\documentclass[12pt, letterpaper]{article}

\usepackage[margin=1in]{geometry}
 
\usepackage{tikz}

\usetikzlibrary{shadings}

\usepackage{amsmath}

\usepackage{amsfonts}

\usepackage{amssymb}

\usepackage{amsthm}

\usepackage[utf8]{inputenc}
\usepackage[T1]{fontenc}
\usepackage{newtxtext}
\usepackage{newtxmath}
\usepackage{bbold}

\usepackage{mathtools}

\usepackage{arydshln}

\usepackage{fancyhdr}

\usepackage{verbatim}

\usepackage{enumitem}

\usepackage[normalem]{ulem}
\usepackage{cancel}

\usepackage{color}

\usepackage[colorlinks=true, linkcolor=teal, citecolor=-teal,
pagebackref=true]{hyperref}

\usepackage{comment}

\usepackage[roundcorner=5pt,  skipabove=10pt, skipbelow=10pt]{mdframed}

\usepackage[colorinlistoftodos]{todonotes}

\newcommand{\thoreau}[1]{\todo[inline, color=cyan!20]{Thoreau says: #1}}
\newcommand{\manuel}[1]{\todo[inline, color=purple!20]{Man{\"u}el says: #1}}
\newcommand{\felipe}[1]{\todo[inline, color=teal!20]{Rambo Ramirez says: #1}}

\newcommand{\demimargin}[1]{\todo[color=orange!20]{Demi says: #1}}
\newcommand{\manuelmargin}[1]{\todo[ color=purple!20]{Manuel says: #1}}
\newcommand{\felipemargin}[1]{\todo[color=teal!20]{Felipe says: #1}}
\usepackage{cleveref}
\usepackage{nameref}

\theoremstyle{plain}
\newtheorem*{theorem*}{Theorem}
\newtheorem{theorem}{Theorem}%[section]
\newtheorem{lemma}{Lemma}[section]

\newtheorem{proposition}[lemma]{Proposition}

\newtheorem*{claim*}{Claim}

\newtheorem*{gallagherthm}{Gallagher's Theorem}

\newtheorem{dhtheorem}{Theorem}
\renewcommand{\thedhtheorem}{DH}
\crefname{dhtheorem}{Theorem}{Theorems}

\newtheorem{bvhtheorem}{Theorem}
\renewcommand{\thebvhtheorem}{BHV}
\crefname{bvhtheorem}{Theorem}{Theorems}

\newtheorem*{celemma}{Chung--Erd{\H o}s Lemma}

\theoremstyle{definition}

\newtheorem*{example*}{Example}

\theoremstyle{remark}
\newtheorem*{remark*}{Remark}

\newcommand{\R}{\mathbb{R}}

\newcommand{\ba}{\mathbf a}

\newcommand{\bA}{\mathbf A}

\newcommand{\balpha}{\boldsymbol \alpha}

\newcommand{\bbeta}{\boldsymbol \beta}

\newcommand{\bb}{\mathbf b}

\newcommand{\bc}{\mathbf c}

\newcommand{\be}{\mathbf e}

\newcommand{\beps}{\boldsymbol \varepsilon}

\newcommand{\bgamma}{\boldsymbol \gamma}
\newcommand{\zq}{\overline{\mathbb{Z}_q^n}}
\newcommand{\zr}{\overline{\mathbb{Z}_r^n}}
\newcommand{\Zq}{\mathbb{Z}_q}

\newcommand{\bn}{\boldsymbol n}

\newcommand{\bp}{\mathbf{p}}

\newcommand{\p}{\mathbf{p}}

\newcommand{\bq}{\mathbf{q}}

\newcommand{\bm}{\boldsymbol m}

\newcommand{\br}{\mathbf r}

\newcommand{\bs}{\mathbf s}

\newcommand{\bt}{\mathbf t}

\newcommand{\bu}{\mathbf{u}}

\newcommand{\bv}{\mathbf{v}}

\newcommand{\bx}{\mathbf x}

\newcommand{\by}{\mathbf y}

\newcommand{\CC}{\mathbb{C}}

\newcommand{\EE}{\mathbb{E}}

\newcommand{\PP}{\mathbb{P}}

\newcommand{\RR}{\mathbb{R}}

\newcommand{\QQ}{\mathbb{Q}}

\newcommand{\NN}{\mathbb{N}}

\newcommand{\N}{\mathbb{N}}

\newcommand{\TT}{\mathbb{T}}

\newcommand{\ZZ}{\mathbb{Z}}

\newcommand{\calA}{\mathcal{A}}

\newcommand{\calB}{\mathcal{B}}

\newcommand{\calC}{\mathcal{C}}

\newcommand{\calD}{\mathcal{D}}

\newcommand{\hatD}{\widehat{D}}

\newcommand{\calF}{\mathcal{F}}

\newcommand{\calG}{\mathcal{G}}

\newcommand{\calH}{\mathcal{H}}

\newcommand{\calL}{\mathcal{L}}

\newcommand{\calM}{\mathcal{M}}

\newcommand{\calP}{P}

\newcommand{\calQ}{\mathcal{Q}}

\newcommand{\calS}{\mathcal{S}}

\newcommand{\calW}{\mathcal{W}}

\newcommand{\va}{\mathbf{a}}
\newcommand{\vb}{\mathbf{b}}
\newcommand{\vc}{\mathbf{c}}

\newcommand{\bone}{\mathbb{1}}

\newcommand{\bzero}{\mathbf{0}}

\renewcommand{\leq}{\leqslant} \renewcommand{\geq}{\geqslant}

\renewcommand{\subset}{\subseteq} \renewcommand{\supset}{\supseteq}

\newcommand{\sprindzuk}{Sprind{\v z}uk}

\DeclarePairedDelimiter{\indcond}{[}{]}
\newcommand{\ind}[1]{\bone_{\indcond*{#1}}}

\DeclarePairedDelimiter{\abs}{\lvert}{\rvert}
\DeclarePairedDelimiter{\inner}{\langle}{\rangle}
\DeclarePairedDelimiter{\norm}{\lVert}{\rVert}
\DeclarePairedDelimiter{\floor}{\lfloor}{\rfloor}
\DeclarePairedDelimiter{\ceil}{\lceil}{\rceil}
\DeclarePairedDelimiter{\set}{\lbrace}{\rbrace}
\DeclarePairedDelimiter{\parens}{\lparen}{\rparen}
\DeclarePairedDelimiter{\brackets}{\lbrack}{\rbrack}

\DeclareMathOperator{\diam}{diam}
\DeclareMathOperator{\denom}{denom}
\DeclareMathOperator{\meas}{m}
\DeclareMathOperator{\outm}{\overline m}
\DeclareMathOperator{\inm}{\underline m}
\DeclareMathOperator{\real}{Re}
\DeclareMathOperator{\imag}{Im}
\DeclareMathOperator{\Log}{Log}
\DeclareMathOperator{\lcm}{lcm}
\DeclareMathOperator{\rep}{rep}

\DeclareMathOperator{\lex}{lex}
\DeclareMathOperator{\supp}{supp}
\DeclareMathOperator{\vol}{vol}

\def\eps{{\varepsilon}}

\def\go{{\omega}}

\def\gs{{\sigma}}

\def\al{{\alpha}}

\def\1int{{[0,1]}}

\def\Sprindzuk{{Sprind{\v z}uk}}

\title{On metric multiplicative Diophantine approximation}

\author{Manuel Hauke-Treuer \\ Graz University of Technology \\ \texttt{hauke@math.tugraz.at} \and Felipe A.~Ram{\'i}rez \\ Wesleyan University \\ \texttt{framirez@wesleyan.edu}}
\date{}

\begin{document}

\maketitle

\begin{abstract}
  Gallagher's theorem (1962) gives a Khintchine-type dichotomy for the
  Lebesgue measure of the set of multiplicatively $\psi$-approximable
  points in $\RR^n$, $n \geq 2$, for monotonic functions $\psi$.  We
  establish this dichotomy for arbitrary functions $\psi$, proving a
  conjecture of Beresnevich, Haynes, and Velani (2013).
\end{abstract}

\section{Introduction}

Multiplicative Diophantine approximation concerns the set
$\calG_n(\psi)$ of
$\balpha = (\alpha_1, \dots, \alpha_n)\in\RR^n$ such that there are infinitely many
solutions to
\begin{equation}\label{eq:multapprox}
  \abs*{\alpha_1 - \frac{a_1}{q}}\cdot  \abs*{\alpha_2 - \frac{a_2}{q}}\cdots   \abs*{\alpha_n - \frac{a_n}{q}} \leq \frac{\psi(q)}{q^n}\qquad (a_1, \dots, a_n, q)\in\ZZ^n\times\NN,
\end{equation}
where $\psi:\NN\to [0, \infty)$ is a given function. The foremost open problem
in this area is Littlewood's conjecture, stating that
\begin{equation*}%\label{eq:littlewood}
  \calL := \bigcap_{\eps>0}\calG_2(\eps/q) = \RR^2.
\end{equation*}
It is a simple matter to show that $\calL$ has full Lebesgue measure, and much
of the work toward Littlewood's conjecture has resulted in further
descriptions of the smallness of the exceptional set
$\RR^2\setminus\calL$. Most notably, Einsiedler, Katok, and
Lindenstrauss~\cite[2006]{EKL06} showed that $\RR^2\setminus\calL$
 is a countable union of sets of zero box-counting
dimension, improving vastly on the baseline knowledge that it is a
null set.

In a related direction of research, one aims to determine the extent
to which the claim in Littlewood's conjecture can be strengthened,
after allowing an exceptional set of Lebesgue measure zero. Gallagher proved
the following theorem.

\begin{gallagherthm}[{\cite[1962]{GallagherMultiplicative}}]
  \phantomsection\label{thm:G}
  Let $n\geq 2$ and let $\psi:\NN\to[0,\infty)$ be monotonically decreasing. Then 
  \begin{equation*}
    \lambda_n(\calG_n(\psi)) =
    \begin{cases}
      0 &\textrm{if } \sum\psi(q)\parens*{\log\frac{1}{\psi(q)}}^{n-1}< \infty, \\[2ex]
      1 &\textrm{if } \sum\psi(q)\parens*{\log\frac{1}{\psi(q)}}^{n-1} = \infty,
    \end{cases}
  \end{equation*}
  where $\log x := \max\{\ln x, 1\}$, and $\lambda_n$ denotes the $n$-dimensional Lebesgue measure on $[0,1]^n$.
\end{gallagherthm}

\begin{remark*}
  In fact, the series appearing in~\cite{GallagherMultiplicative} is
  $\sum\psi(q)(\log q)^{n-1}$, but the above statement is equivalent,
  in the presence of the monotonicity assumption, see
    e.g. \cite[Remark 4]{BHV}. We remark that the above theorem is
  also true for $n=1$. In this case it is exactly Khintchine's
  theorem~\cite[1924]{Khintchineonedimensional}, and hence, Gallagher's theorem
  is the extension of Khintchine's theorem to the multiplicative setting.
\end{remark*}

Gallagher's theorem implies that Littlewood's conjecture can be
strengthened by a factor of $(\log q)^2$ for almost all points, but
for almost no points by a factor of
$(\log q)^{2+\delta}$, that is,
\begin{equation*}
  \lambda_2\parens*{\bigcap_{\eps>0} \calG_2\parens*{\frac{\eps}{q(\log q)^2}}} = 1,
  \quad
  \textrm{while}
  \quad
  \lambda_2\parens*{\bigcap_{\eps>0} \calG_2\parens*{\frac{\eps}{q(\log q)^{2+\delta}}}} = 0,
\end{equation*}
for every $\delta>0$. Similarly, it implies that the $n$-dimensional
analogue of Littlewood's conjecture, which is open for every $n \geq 2$, holds almost
everywhere in $\RR^n$ with an extra factor of $(\log q)^n$, but almost
nowhere with an extra factor of $(\log q)^{n+\delta}$. Of course,
Gallagher's theorem says much more, giving a complete metric dichotomy
for monotonic rates, in the style of Khintchine's
theorem~\cite[1924]{Khintchineonedimensional}. Famously, the
monotonicity assumption cannot be removed from Khintchine's
theorem~\cite[1941]{duffinschaeffer}. In contrast, it has been conjectured by
Beresnevich, Haynes, and Velani that Gallagher's theorem holds for
\emph{arbitrary} functions~\cite[2013, Conjecture~1]{BHV}. In this paper, we prove
that conjecture.
  
\begin{theorem}\label{main_thm}
Let $n\geq 2$ and let $\psi:\NN\to[0,\infty)$ be arbitrary. Then 
  \begin{equation*}
    \lambda_n(\calG_n (\psi)) = 
    \begin{cases}
      0 &\textrm{if } \sum\psi(q)\parens*{\log \frac{1}{\psi(q)}}^{n-1}  < \infty,  \\[2ex]
      1 &\textrm{if } \sum\psi(q)\parens*{\log \frac{1}{\psi(q)}}^{n-1}  = \infty.
    \end{cases}
  \end{equation*}
\end{theorem}

Similar questions have been considered for other forms of
approximation, such as simultaneous, dual, systems of linear forms, as
well as inhomogeneous approximation and approximation with coprimality
restrictions~\cite{allen2026inhomogeneouskhintchinetheoremdimension,allen2021independence,allen2023inhomogeneous,BVKG,chow2023littlewood,duffinschaeffer,FH,Gallagherkt,Harman,HaukeRamirezDiscreteAnalysis,KM,PollingtonVaughan,ramirez2023duffinschaeffer,Sprindzuk}. The
case of \textit{simultaneous} approximation was handled by Gallagher in a
follow-up paper~\cite{Gallagherkt} to the above
cited~\cite{GallagherMultiplicative}. He showed that the
higher-dimensional Khintchine theorem~\cite[1926]{Khintchine} can be
stated with arbitrary functions.

The main difficulty in the \emph{multiplicative} setup comes from the
unwieldy geometry of the sets involved. Specifically, the
multiplicative problems present nonconvex sets whose intersections one
must
understand. (See~\Cref{fig:AqHq,fig:Aqprime,fig:Aqtilde,fig:errorcontrib}.)
By contrast, the other forms of approximation mentioned above present
metric balls and thickened subspaces, whose intersections are more
tractable.  Historically, the literature has exhibited a trade-off
between convexity and monotonicity. Beresnevich, Haynes, and Velani
give an enlightening discussion of this trade-off and various
instances of it:
\begin{quote}
  \textit{``[T]he current body of metrical results for lim sup sets
    requires that either the approximating domains are convex or that
    the approximating function is monotonic.''}

  \hfill\cite[Section~1.1]{BHV}
\end{quote}
Khintchine's theorem for simultaneous approximation in higher
dimensions~\cite[1926]{Khintchine} used both convexity and
monotonicity. Later, Gallagher~\cite{Gallagherkt} showed that
monotonicity could be removed from that theorem; his arguments relied
heavily on the convexity of the sets that the simultaneous
approximation setup provides. By contrast, the proof of the multiplicative, and hence non-convex, \hyperref[thm:G]{Gallagher
    theorem}~\cite{GallagherMultiplicative} uses monotonicity in a
  crucial way.

Beresnevich, Haynes, and Velani also point out that what is commonly
called Gallagher's theorem is actually a consequence of a more general
statement established in~\cite{GallagherMultiplicative}:
  \begin{quote}
   \textit{``In fact, Gallagher considers lim sup sets arising from more
    general domains but monotonicity plays a crucial role in his
    approach and seems to be unremovable.''}

  \hfill\cite[Section~1.1]{BHV}
  \end{quote}
Gallagher proved in~\cite{GallagherMultiplicative} that instead of assuming convexity, a Khintchine-type behaviour 
under monotonicity holds for $q\boldsymbol{\alpha} \pmod 1 \in S_q$
 for infinitely many $q$, as long as the sets $(S_q)_{q \in \N}$ satisfy
\begin{equation}\label{eq:starproperty_}
  \bx=(x_1, \dots, x_n) \in S_q \quad\textrm{and}\quad |x_i'|\leq |x_i|\quad  (i=1, \dots, n) \quad \implies \quad  \bx'=(x_1', \dots, x_n') \in S_q.
\end{equation}
The multiplicative approximation sets $q\overline{H}_q$
(see~\Cref{fig:AqHq}) satisfy \eqref{eq:starproperty_}, and
thus~\hyperref[thm:G]{Gallagher's multiplicative theorem} is indeed an
instance of this more general statement. A straightforward
construction confirms the remark quoted above by showing that if one
only assumes~(\ref{eq:starproperty_}), then monotonicity is
unremovable.

\begin{theorem}\label{thm:nostar_}
  For every $n\geq 1$, there exists a sequence of sets $S_q$ 
  satisfying~(\ref{eq:starproperty_}) such that $\sum\lambda(S_q)$ diverges and
  such that the set $\calW$ of points $\balpha\in[0,1]^n$ for which there are infinitely
  many solutions to
  \begin{equation*}
    q\balpha\mod 1 \in S_q, \qquad q\in\NN
  \end{equation*}
  has Lebesgue measure $0$. 
\end{theorem}

In fact, the sets we construct for the proof of Theorem
\ref{thm:nostar_} are even convex. From a geometric perspective,
removing monotonicity might require the consideration of further
properties of the sets $S_q$, perhaps a bound on the ratio between the
diameter (of reasonably truncated subsets of $S_q$) and the volume:
Theorem~\ref{thm:nostar_} is based on very thin rectangles that reduce
the problem to the one-dimensional setting, where monotonicity is
known to be necessary. As we will see in the proof of Theorem \ref{main_thm}, the hyperbolic shapes for multiplicative
approximation are \textit{just} ``higher-dimensional enough'' to allow
the removal of monotonicity.  The upshot is that the multiplicative
problem~(\ref{eq:multapprox}) strikes a rare balance that allows us to
dispense with both convexity and monotonicity for~\Cref{main_thm}.

\subsection*{Organization of the paper}
The remainder of the paper is organized as follows. In Section
\ref{sec:proof-strategy} we present the proof strategy in the context of previous methods, highlighting how our proof differs from previous approaches.  In~\Cref{sec:notation} we
collect notation, and in Section \ref{sec:tools-from-liter} we state
the established tools from the literature that are used in the 
proof of Theorem~\ref{main_thm}. In Section \ref{prelim_lem} we make
elementary calculations that are directly used in the subsequent technical sections. In Section \ref{sec:overlap} we provide the
technical heart of the paper by calculating upper bounds on second
moments, and in Section \ref{sec:final} we conclude the proof of
Theorem~\ref{main_thm}, as well as the (particularly
  short) proof of Theorem~\ref{thm:nostar_}.

\subsection*{Acknowledgments}
MHT was funded by the Austrian Science Fund (FWF) project
10.55776/ESP5134624. Parts of this project were carried out while MHT
visited FAR. MHT would like to thank Wesleyan University, as well as
the whole Ram{\'i}rez family (Yeni, Luna, Joaqu{\'i}n and Felipe) for
their hospitality. More precisely, he thanks Yeni for the excellent
Sancocho de Cola, Luna for teaching him how to do Origami, and
Joaqu{\'i}n and Felipe for showing him where to find Asteroids, Waldo,
and other things that really matter in life.

\section{The proof strategy}\label{sec:proof-strategy}

\begin{figure}[t]
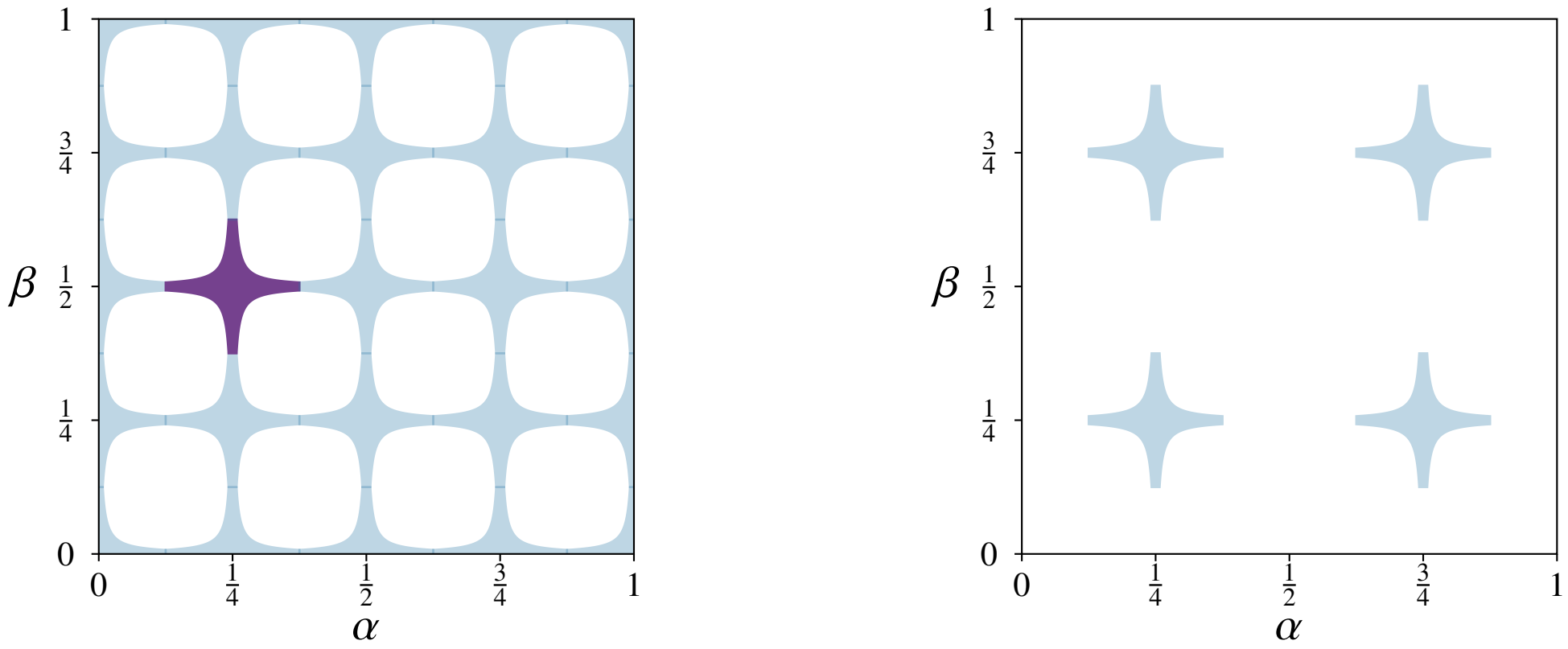

  \centering

  \caption{On the left: An illustration of a set $A_q$ with one of the translated copies of $\overline{H}_q$ highlighted in purple. In this example,
    $q=4$, and the highlighted set is
    $\overline{H}_4+(\frac14,\frac12)$. On the right: An illustration
    of the sets $A_q''$ appearing in
    Gallagher~\cite{GallagherMultiplicative}, coming from the
    condition that each coordinate of $\ba$ is coprime to
    $q$. Gallagher proved QIA for these sets, under a monotonicity
    assumption.}
  \label{fig:AqHq}
\end{figure}

In this section, we outline the strategy of the proof of Theorem \ref{main_thm}, focusing on the case
where $n = 2$.  Note that $\calG_2 (\psi) = \limsup_{q\to \infty}A_q$ where
\[
  A_q = \dot{\bigcup_{\mathbf{a} \in \mathbb{Z}_q^2}} \overline{H}_q + \frac{\mathbf{a}}{q},
\]
with
\[
  \overline{H}_q + \frac{\mathbf{a}}{q} = \left\{(\alpha,\beta) \in
    [0,1]^2: \left\lvert \alpha - \frac{a_1}{q}\right\rvert\cdot\left\lvert
      \beta - \frac{a_2}{q}\right\rvert \leq \frac{\psi(q)}{q^2},
    \left\lvert \alpha - \frac{a_1}{q}\right\rvert \leq
    \frac{1}{{2}q}, \left\lvert \beta - \frac{a_2}{q}\right\rvert \leq
    \frac{1}{{2}q}\right\}.
\]
See~\Cref{fig:AqHq} for an illustration of $A_q$ and $\overline{H}_q$.
\suppressfloats[t]

A straightforward calculation (Lemma~\ref{lem:Aqmeasure}) shows that
$\lambda(A_q) \asymp \psi(q)\log \frac{1}{\psi(q)}$ and, therefore,
the convergence part of~\Cref{main_thm} follows by an application of
the first Borel--Cantelli lemma. The divergence part
of~\Cref{main_thm} would follow from the second Borel--Cantelli lemma,
if the sets $(A_q)_{q \in \mathbb{N}}$ were stochastically
independent. However, they are far from independent.  Instead, we aim
to show that they contain subsets
$A_q^*\subset A_q$ such that $\sum\lambda(A_q^*)$
  diverges and such that the $(A_q^*)_q$ exhibit a weaker form of
independence, called quasi-independence on average (QIA).
Specifically, it would suffice to prove that there exists an absolute
constant $C > 0$, and subsets $A_q^* \subseteq A_q$ such that for
infinitely many disjoint ranges $[X,Y]$ with
$\sum_{X\leq q \leq Y} \lambda(A_q^*) \in [1,2]$, we have
\begin{equation}\label{QIA_model}
  \sum_{q,r \in [X,Y]} \lambda(A_q^* \cap A_r^*) \leq C \sum_{q,r \in [X,Y]} \lambda(A_q^*)\lambda(A_r^*) \leq 4C,
\end{equation}
for then the~\hyperref[thm:CE]{Chung--Erd{\H o}s Lemma}
(see~\Cref{sec:tools-from-liter}) would imply that in these ranges we
had
\begin{equation*}
  \lambda\left(\bigcup_{X\leq q\leq Y} A_q\right) \geq \lambda\left(\bigcup_{X\leq q\leq Y} A_q^*\right) \geq (4C)^{-1},
\end{equation*}
and $\lambda(\calG_2(\psi))\geq (4C)^{-1}$ would follow.  A zero-one law
(\Cref{thm:BHV}, see~\Cref{sec:tools-from-liter}
and~\cite[Theorem~1]{BHV}) would then give \Cref{main_thm}.

The most straightforward approach when aiming for \eqref{QIA_model}
would be to choose $A_q^* = A_q$.  Computing $\lambda(A_q \cap A_r)$ amounts to
summing the volume of
$(\overline{H}_q + \frac{\bf{a}}{q}) \cap (\overline{H}_r +
\frac{\bf{b}}{r})$ over all values of
$(\mathbf{a},\mathbf{b}) \in \mathbb{Z}_q^2 \times \mathbb{Z}_r^2$. Clearly, the volume of
$(\overline{H}_q + \frac{\bf{a}}{q}) \cap (\overline{H}_r +
\frac{\bf{b}}{r})$ is determined by
$\frac{\bf{a}}{q} - \frac{\bf{b}}{r}$. Hence, writing
\begin{equation*}
  \hat{w}(\mathbf{x}) = \hat{w}_{q,r}(\mathbf{x}):= \lambda\left((\overline{H}_q + \mathbf{x}) \cap \overline{H}_r\right),
\end{equation*}
one has
\begin{equation}\label{eq:1}
  \lambda(A_q \cap A_r)= \sum_{(\mathbf{a},\mathbf{b}) \in \mathbb{Z}_q^2 \times \mathbb{Z}_r^2}
  \hat{w}\left(\frac{\bf{a}}{q} - \frac{\bf{b}}{r}\right) =
  \gcd(q,r)^2 \sum_{\mathbf{c} \in \mathbb{Z}_{\lcm(q,r)}^2}
  \hat{w}\left(\frac{\mathbf{c}}{\lcm(q,r)}\right),
\end{equation}
by an application of the Chinese Remainder Theorem. Since
$\hat w(\bx)$ is monotonically decreasing in $|x_i|$ for $i = 1,2$ (here, the property \eqref{eq:starproperty_} is decisive), we may bound the sum with
an integral,
\begin{equation}\label{eq:2}
  \sum_{\bc \in \ZZ^2}\hat{w}\parens*{\frac{\bc}{\lcm(q,r)}} \leq \lcm(q,r)^2 \int_{\RR^2} \hat{w}(\bx)\, \mathrm{d}\bx + \sum_{\substack{\bc \in \ZZ^2 \\ c_1 = 0 \textrm{ or } c_2=0}}\hat{w}\parens*{\frac{\bc}{\lcm(q,r)}}.
\end{equation}
Straightforward calculations reveal that
\begin{equation}\label{eq:3}
  \int_{\mathbb{R}^2}\hat{w}\left(\mathbf{x}\right)\mathrm{d}\bx \asymp \frac{1}{q^2}\frac{1}{r^2} \lambda(A_q) \lambda(A_r),
\end{equation}
which, combined with~(\ref{eq:1}),~(\ref{eq:2}), and~(\ref{eq:3}), shows
\begin{equation}\label{eq:4}
  \lambda(A_q\cap A_r) \ll \lambda(A_q)\lambda(A_r) + \gcd(q,r)^2 \sum_{\substack{\bc \in \ZZ^2 \\ c_1 = 0 \textrm{ or } c_2=0}}\hat{w}\parens*{\frac{\bc}{\lcm(q,r)}}.
\end{equation}
Evidently, the only trouble can come from pairs $(\ba, \bb)$ where
\begin{equation}\label{eq:coordcoincidence}
  \frac{a_1}{q} = \frac{b_1}{r}\qquad\textrm{or} \qquad   \frac{a_2}{q} = \frac{b_2}{r}.
\end{equation}
However, calculations show that the contribution from these
points can indeed destroy the desired QIA property~(\ref{QIA_model}) for the choice $A_q^* = A_q$.

Gallagher's idea~\cite{GallagherMultiplicative} was to show
\eqref{QIA_model} for $A_q'' \subseteq A_q$ instead, defined by
\[
  A_q'' := \bigcup_{\substack{\textbf{a} = (a_1,a_2) \in \mathbb{Z}_q^2,\\
      \gcd(a_1,q) = 1,\\
      \gcd(a_2,q) = 1}} \left(\overline{H}_q + \frac{\mathbf{a}}{q}\right).
\]
(See~\Cref{fig:AqHq}.) For $q\neq r$, the coprimality conditions rule out the
existence of contributions from pairs $(\ba,\bb)$
satisfying~(\ref{eq:coordcoincidence}). However, one now has
$\lambda(A_q'') \asymp
\left(\frac{\varphi(q)}{q}\right)^2\lambda(A_q)$ (where $\varphi$ denotes the Euler totient function), and if $\psi$ is
supported only on numbers $q$ with many small prime factors, this can
cause the series $\sum\lambda(A_q'')$ to converge. It is for this
reason that Gallagher's theorem has a monotonicity assumption on
$\psi$: the weights $\varphi(q)/q$ are
bounded away from $0$ on average, so employing partial
    summation shows that for monotonic $\psi$ we have
  \begin{equation}
    \label{eq:quickcomp}
    \sum_{q\leq Q}\lambda(A_q'')\asymp \sum_{q\leq Q}\lambda(A_q),
  \end{equation}
  and the series $\sum\lambda(A_q'')$ diverges.
Now, having removed the contributions from~(\ref{eq:coordcoincidence})
 the bound~(\ref{eq:4}) becomes
\[
  \lambda(A_q'' \cap A_r'') \ll \lambda(A_q)\lambda(A_r), 
\]
and~(\ref{eq:quickcomp}) is used again to show that one has
\[
  \lambda(A_q)\lambda(A_r) \ll \lambda(A_q'')\lambda(A_r''),
\]
on average over long ranges of $q,r$. This
  implies~\eqref{QIA_model} which, as discussed, is sufficient for
  full measure of the limsup set. This proves \hyperref[thm:G]{Gallagher's
  theorem}, but having crucially used the monotonicity assumption. 

\suppressfloats[t]
\begin{figure}[t]
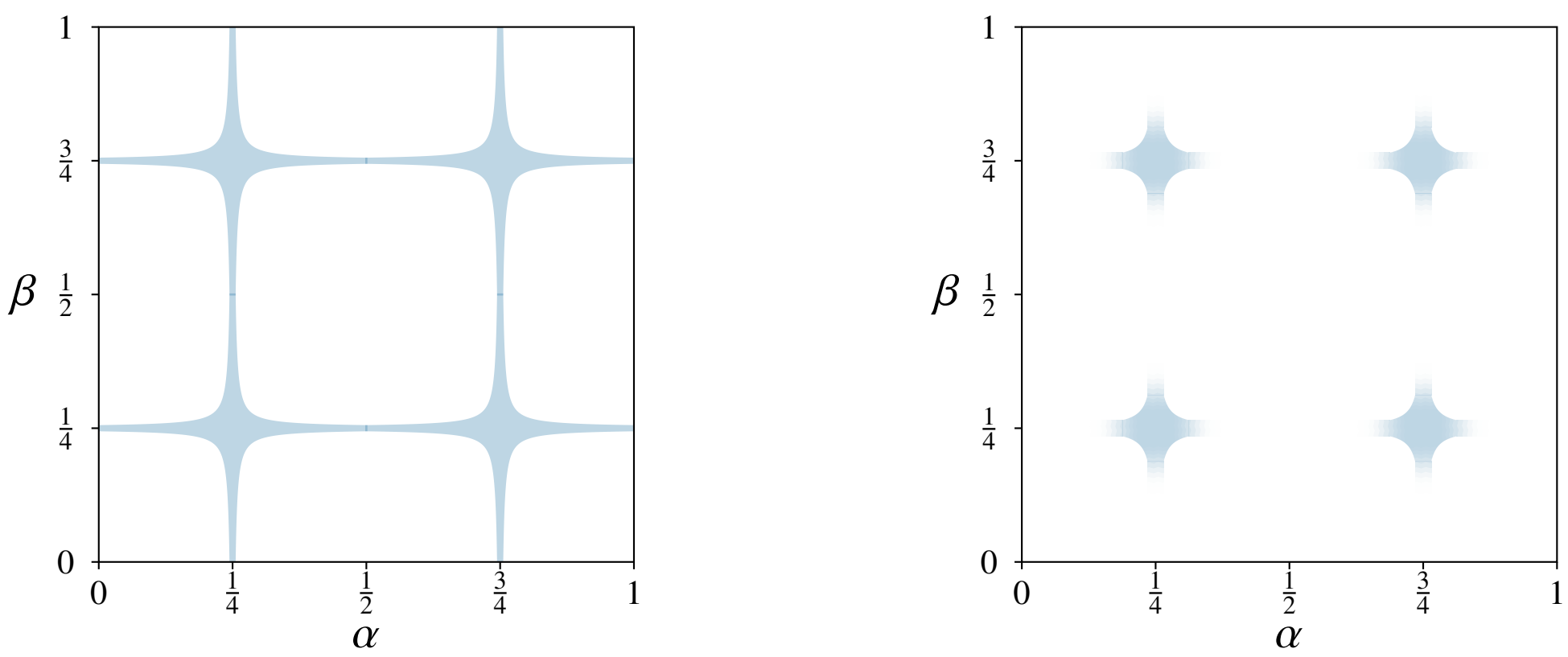

  \centering

  \caption{The sets $A_q'$ and $A_q^{*,\prime}$ appearing in
    Fr{\"u}hwirth--Hauke~\cite{FH}. Note the long arms on $A_q'$ in
    comparison to the short arms of $A_q''$ in~\Cref{fig:AqHq}, and
    truncations in $A_q^{*,\prime}$.}
  \label{fig:Aqprime}
\end{figure}

In fact, it was shown by Fr{\"u}hwirth and the first author~\cite{FH}
that a coordinate-wise coprimality condition in~(\ref{eq:multapprox})
also leads to Khintchine-type behaviour, without assuming monotonicity: denoting the underlying sets by $A_q'$, it was shown that
$\sum_{q \in \mathbb{N}}\lambda(A_q') = \infty$ implies
$\lambda(\limsup_{q \to \infty} A_q') = 1$ for arbitrary $\psi$. This resolved a
conjecture of Beresnevich, Haynes, and
Velani~\cite[Conjecture~2]{BHV}, and incidentally clarified an
overlooked detail in the problem, namely, that $A_q'$ has ``arms'' that
extend past the arms of $A_q''$ (see~\Cref{fig:Aqprime}), and result in
a different divergence condition from the conjectured one appearing in~\cite{BHV}.
As in the previously described argument for Gallagher's theorem, it
was necessary to restrict to subsets of the sets arising naturally for the respective approximation sets. Instead of working directly with $A_q'$,
 the arguments in~\cite{FH} used fiberwise-smooth bump
  functions $0\leq f_q \leq \bone_{A_q'}$ that could heuristically be thought of
  as producing truncated sets $A_q^{*,\prime}\subset A_q'$ similar to the
  ones~\Cref{fig:Aqprime} depicts.\footnote{These in turn facilitated
    a fiberwise application of the Koukoulopoulos--Maynard Theorem. For
    details, we refer readers to~\cite{FH}.}

\begin{figure}[t]
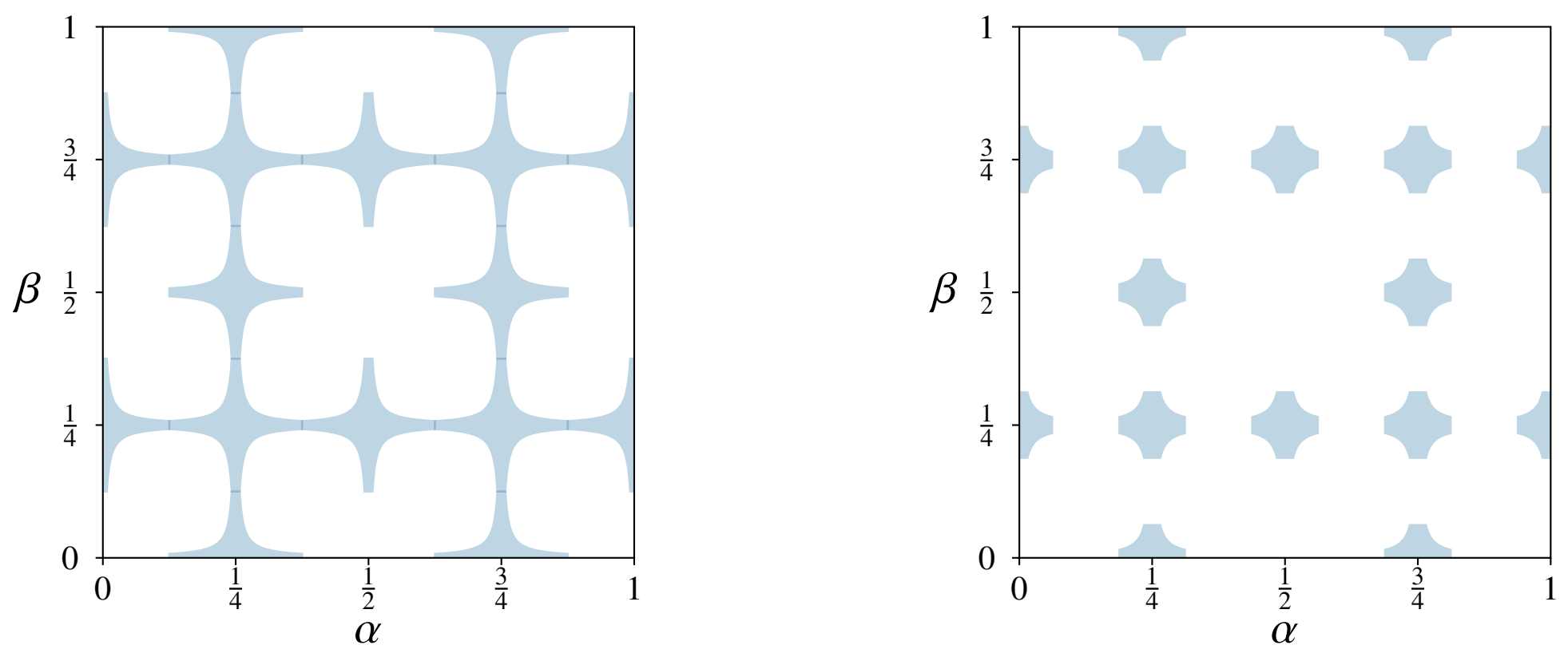

  \centering

  \caption{On the left: The sets $\tilde{A}_q$ defined using
    $\gcd(a_1,a_2,q)=1$. On the right: Their truncations $A_q^*$.}
  \label{fig:Aqtilde}
\end{figure}

Returning to the original question, which assumes only the divergence
of $\sum_{q}\lambda(A_q)$, it is clear that restricting to $A_q''$
will not help, because without monotonicity, there is no guarantee that
the series $\sum\lambda(A_q'')$ diverges. We will see, however, that
instead of \textit{pairwise} coprimality
$\gcd(a_1,q) = \gcd(a_2,q) = 1$, we can use \textit{jointly} coprime
vectors, i.e.,
\[
  \tilde{A}_q := \bigcup_{\substack{\textbf{a} = (a_1,a_2) \in \mathbb{Z}_q^2,\\
\gcd(a_1,a_2,q) = 1,\\
}} \left(\overline{H}_q + \frac{\mathbf{a}}{q}\right).
\]
(See~\Cref{fig:Aqtilde}.) Unlike coprimality in dimension one or the pairwise coprimality assumption employed by Gallagher,
  the joint coprimality condition $\gcd(a_1,a_2,q) = 1$ holds for a positive proportion of
  $\textbf{a} \in \mathbb{Z}_q^2$ for \textit{every} $q$. Therefore, the divergence of
  $\sum_q\lambda(A_q)$ implies the divergence of
  $\sum_q\lambda(\tilde{A}_q)$ without having to assume monotonicity.
For $q \neq r$, the condition $\gcd(a_1,a_2,q) = 1$ rules out coincidences where
\emph{both} equalities in~(\ref{eq:coordcoincidence}) hold,  that is,
\begin{equation*}
  \frac{a_1}{q} = \frac{b_1}{r}\qquad\textrm{AND} \qquad   \frac{a_2}{q} = \frac{b_2}{r}, \tag{i}
\end{equation*}
but cannot avoid contributions from coordinate-wise coincidences
\begin{equation*}
  \frac{a_1}{q} = \frac{b_1}{r}\; \textrm{ and }\; \frac{a_2}{q} \neq \frac{b_2}{r}
  \qquad\textrm{OR}\qquad
\frac{a_1}{q} \neq \frac{b_1}{r}\; \textrm{ and }\; \frac{a_2}{q} = \frac{b_2}{r}.
   \tag{ii} 
\end{equation*}
Avoiding (i) (but not (ii)), allows us to refine~(\ref{eq:2}) by bounding
\begin{equation}\label{one_dim_int}
\sum_{\substack{\bc \in \ZZ^2\setminus\set{\bzero} \\ c_1 = 0 \textrm{ or } c_2=0}}\hat{w}\parens*{\frac{\bc}{\lcm(q,r)}} \ll \lcm(q,r) \int_{\RR} \hat{w}(x,0)\,\mathrm{d}x,
\end{equation}
again using monotonicity and symmetry of $\hat w$. The
arguments that led to~(\ref{eq:4}) now give
\begin{align*}
  \lambda(\tilde{A}_q\cap \tilde{A}_r) &\ll \lambda(A_r)\lambda(A_q) + \gcd(q,r)^2\lcm(q,r)  \int_{\RR} \hat{w}(x,0)\,\mathrm{d}x \\
              &\ll \lambda(\tilde{A}_q)\lambda(\tilde{A}_r) + \underbrace{qr\gcd(q,r)  \int_{\RR} \hat{w}(x,0)\,\mathrm{d}x}_{:= \hat{E}_{q,r}}.
\end{align*}
Evidently, we must study the error terms $\hat{E}_{q,r}$. If one could
show that
\begin{equation}\label{tilde_e_model}
  \sum_{q,r \in [X,Y]}\hat{E}_{q,r} \ll \left(\sum_{q \in [X,Y]}\lambda(A_q)\right)^2 \asymp 1,
\end{equation}
then the sufficient QIA \eqref{QIA_model} for $(\tilde{A}_q)_q$ would follow. 

We note that Gallagher himself used the same joint primitivity
condition to remove monotonicity from Khintchine's
Theorem~\cite{Khintchine} in the question of \textit{simultaneous}
 approximation~\cite{Gallagherkt}. There, the underlying approximation
sets are squares, hence the function $w$ behaves much more
regularly. In fact, a direct application of the above method with
$\hat{w}(x,0) \leq \lVert\hat{w}\rVert _{\infty}$ allows us to
establish \eqref{tilde_e_model} for the analogous sets. This highlights the convenience of working with sets having nice
  geometric properties. Khintchine's theorem is about balls in a
  norm-induced metric; they are convex and have bounded
  eccentricity.

For the actual sets $(\tilde{A}_q)_q$ and the underlying hyperbolic
regions $\overline{H}_q$, this crude approach fails. In order to study
\eqref{tilde_e_model} in more detail in an idealized setup, we pose
the following Model problem: Consider a set
$\mathcal{M} \subseteq [X,2X]$ with $\#\mathcal{M} = N$, where $X$ is
arbitrarily large in terms of $N$ (in other words, assume $Y = 2X$ in
the above). Further, assume that $\supp \psi = \mathcal{M}$, $\psi$ is
constant on its support, and the total measure is comparable to $1$,
i.e.,
\[
  \sum_{q \in \mathcal{M}} \psi(q) \log\frac{1}{\psi(q)} = 1,
\]
which implies
\[
  \psi(q) \asymp \frac{1}{N \log N},\quad q \in \mathcal{M}.
\]
The Model problem is to show~(\ref{tilde_e_model}) in this simplified scenario.
Now it is particularly easy to compute that
\[
\hat{E}_{q,r}\asymp  \gcd(q,r)\min\left\{\frac{\psi(q)}{q},\frac{\psi(r)}{r}\right\}
  \asymp \frac{1}{N \log N} \frac{\gcd(q,r)}{\sqrt{qr}}.
\]
However, there are known ensembles $\mathcal{M} \subset \N$ of cardinality $N$ where
\[\sum_{q,r \in \mathcal{M}} \frac{\gcd(q,r)}{\sqrt{qr}} = \omega(N \log N),\]
thus~(\ref{tilde_e_model}) fails. Note that this is not due to
coarseness of estimates, but a second moment that is really too
large.

\begin{figure}[h!]
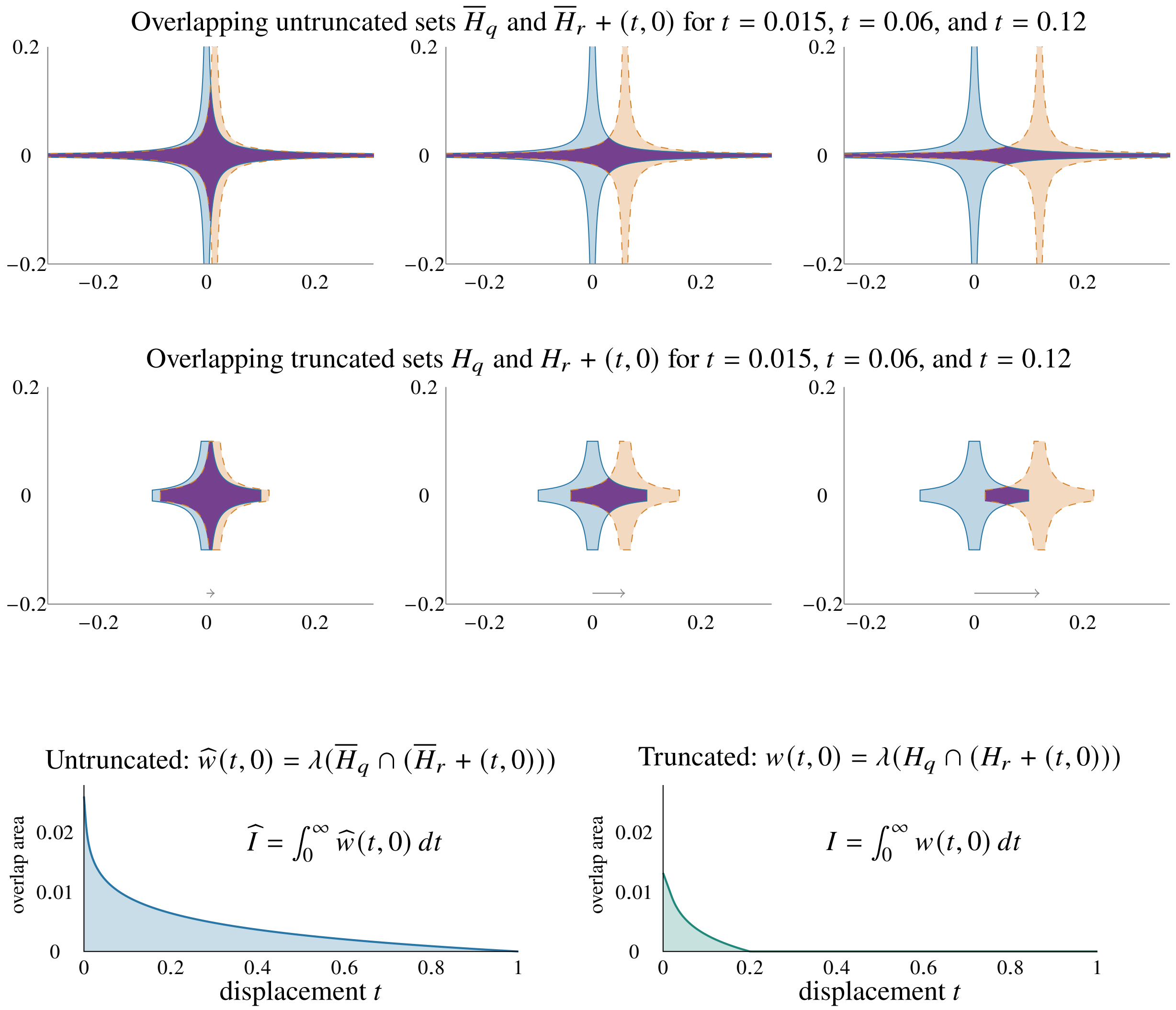

  \centering

  \caption{The contribution of the error terms. The top row shows the
    overlaps between an untruncated set $\overline{H}_q$ and
    coordinatewise shifts of $\overline{H}_r$ when there
      is a coincidence in the other coordinate. The middle row shows
    the same for the truncated sets $H_q$ and $H_r$. Note the smaller
    purple regions of intersection for the truncated sets. The bottom
    row shows a comparison of the integrals of the measures of those
    overlaps.}
  \label{fig:errorcontrib}
\end{figure}

In order to overcome the issue, we further prune the sets
$(\tilde{A}_q)_q$, but this time at a local level: for fixed
$\textbf{a} = (a_1,a_2)$, instead of considering the whole hyperbolic
region
\[
  \overline{H}_q + \frac{\mathbf{a}}{q} =  \left\{(\alpha,\beta) \in \mathbb{R}^2: \left\lvert \alpha - \frac{a_1}{q}\right\rvert\cdot \left\lvert \beta - \frac{a_2}{q}\right\rvert \leq \frac{\psi(q)}{q^2}, \left\lvert \alpha - \frac{a_1}{q}\right\rvert \leq \frac{1}{2q}, \left\lvert \beta - \frac{a_2}{q}\right\rvert \leq \frac{1}{2q}\right\},
\]
we consider the \textit{truncated} region (see~\Cref{fig:Aqtilde})
\[
  H_q + \frac{\mathbf{a}}{q}  =  \left\{(\alpha,\beta) \in \mathbb{R}^2: \left\lvert \alpha - \frac{a_1}{q}\right\rvert\cdot \left\lvert \beta - \frac{a_2}{q}\right\rvert \leq \frac{\psi(q)}{q^2}, \left\lvert \alpha - \frac{a_1}{q}\right\rvert \leq \frac{\psi(q)^{1/3}}{{2}q}, \left\lvert \beta - \frac{a_2}{q}\right\rvert \leq \frac{\psi(q)^{1/3}}{{2}q}\right\}.
\]
Note that this truncation satisfies $\lambda(H_q) \asymp \lambda(\overline{H}_q)$ (see Lemma \ref{lem:Aqmeasure}), which consequently implies that for
\[A_q^* := \bigcup_{\substack{\mathbf{a} \in \mathbb{Z}_q^2\\ \gcd(a_1,a_2,q) = 1}}\left( H_q + \frac{\mathbf{a}}{q}\right),\]
 we still have 
\[\lambda(A_q^*) \asymp \lambda(\tilde{A}_q) \asymp \lambda(A_q).\]
  This pruning results in a new function
\begin{equation*}
  w(\mathbf{x}) = w_{q,r}(\mathbf{x}):= \lambda\left((H_q + \mathbf{x}) \cap H_r\right),
\end{equation*}
that has substantially smaller values along the coordinate axes than
the original $\hat{w}$ did. (See~\Cref{fig:errorcontrib}.)
The error term to control is now
\begin{equation*}
  E_{q,r} := qr \gcd(q,r) \int_{0}^\infty w(0,x) \,\mathrm{d}x.
\end{equation*}
A slightly longer yet elementary calculation (Lemma~\ref{lem:intHqHr}) shows that
\[
  E_{q,r} \ll \frac{\gcd(q,r)}{\sqrt{qr}}\min\left\{\psi(q)^{7/6},
    {\psi(r)^{7/6}} \right\} \asymp \frac{1}{(N \log N)^{7/6}}
  \frac{\gcd(q,r)}{\sqrt{qr}}.
\]
By a classical GCD sum bound, see~\Cref{lem:gcdsum},
\[\sum_{q,r \in \mathcal{M}} \frac{\gcd(q,r)}{\sqrt{qr}} \ll N^{1+ \varepsilon},\]
which implies that 
\[
  \sum_{q,r \in \mathcal{M}} E_{q,r} \ll \frac{1}{(N \log N)^{7/6}} \sum_{q,r \in
    \mathcal{M}}\frac{\gcd(q,r)}{\sqrt{qr}} \ll \frac{N^\varepsilon}{N^{1/6}} \ll 1.
\]
This established \eqref{QIA_model} in the Model problem, where the support was restricted to dyadic ranges with constant function values on its support.

Going from the Model problem to arbitrary sets $\mathcal{M}$ with
varying values of $\psi$ is mainly a technical challenge that must be
undertaken in order to correctly evaluate, or at least efficiently bound, an error
term $E_{q,r}$ that now depends non-trivially on $q$ and $r$ as well as
$\psi(q)$ and $\psi(r)$. However, an application of Cauchy--Schwarz reduces this situation essentially to the Model problem. The passage to the higher-dimensional
case $n \geq 3$ also works in a reasonably straightforward manner: instead of primitive vectors, we define the sets $\tilde{A}_q$ using vectors
$(a_1,\ldots,a_n)$ where $\gcd(a_i,a_j,q) = 1$ for all $i \neq j$, a
condition which is again satisfied by a positive proportion of all
vectors $\mathbf{a} \in \mathbb{Z}_q^n$. Now the $n$-dimensional version
of the error term $\hat{E}_{q,r}$ contains an $(n-1)$-dimensional
integral $\int_{\mathbb{R}^{n-1}} \hat{w}(0,\mathbf{x}) \mathrm{d}\bx$
whose value is too large for QIA and, again, the issue is fixed by
restricting to polynomially truncated sets $A_q^*\subset
\tilde{A}_q$. (See~\Cref{lem:intHqHr}.)

\section{Notation}\label{sec:notation}

Throughout, $\psi:\NN\to [0,\infty)$ is an arbitrary function. The quantity
$\psi(q))\parens*{\log\frac{1}{\psi(q)}}^{n-1}$ is of special interest, so
we reiterate the conventions that $\log x = \max\{1,\ln x\}$ and
$0\cdot\log(1/0)^{n-1} = 0$.

We study the inequality
\begin{equation}\label{eq:inequality}
  \abs*{\alpha_1 - \frac{a_1}{q}}\cdot  \abs*{\alpha_2 - \frac{a_2}{q}}\cdots   \abs*{\alpha_n - \frac{a_n}{q}} \leq \frac{\psi(q)}{q^n}\qquad (a_1, \dots, a_n, q)\in\ZZ^n\times\NN,
\end{equation}
where $\balpha = (\alpha_1,\dots, \alpha_n) \in\RR^n$ and, in particular, the set
\begin{equation*}
    \calG_n (\psi) = \set*{\balpha\in [0,1]^n : \textrm{\eqref{eq:inequality} holds for infinitely many $(\ba,q)\in\ZZ^n\times \NN$}}. 
\end{equation*}
Note that we now adopt the convention of viewing $\calG_n$ as a subset of $[0,1]^n$. No generality is lost, since its defining condition is $\ZZ^n$-periodic.
We define
\begin{align*}
  G_q &:=G_q^{(n)}(\psi) = \set*{(x_1,\dots, x_n)\in \brackets*{-\frac{1}{2}, \frac{1}{2}}^n: \prod_{i}|x_i| \leq \psi(q)},\\
  H_q &:=H_q^{(n)}(\psi) = \frac{1}{q}\parens*{G_q^{(n)}(\psi) \cap \brackets*{-\frac{\psi(q)^{\frac{1}{n+1}}}{2}, \frac{\psi(q)^{\frac{1}{n+1}}}{2}}^n}, \\
  A_q &:=A_q^{(n)}(\psi) = \{\balpha\in[0,1]^n: (\ref{eq:inequality}) \textrm{ holds  for some } \ba \in\ZZ^n \}, \\
  A_q^* &:=A_q^{*,(n)}(\psi) = \bigcup_{\substack{\ba\in\ZZ_q^n \\ \gcd(a_i,a_j,q)=1\, (i\neq j)}} \left(\frac{\ba}{q} + H_q\right),
\end{align*}
and note that $\limsup A_q^*\subset \limsup A_q = \calG_n(\psi)$.

We use $\lambda_n$ to denote Lebesgue measure on $[0,1]^n$, often
dropping the subscript when the context makes it clear, and write
$\mathbb{Z}_n := \mathbb{Z}/n\mathbb{Z}$.  We use the classical
$O$-notation as well as the Vinogradov notation $\ll$ for
$f \ll g \Leftrightarrow f = O(g) \Leftrightarrow \exists c> 0: f(x)
\leq c g(x)$, and write $f\asymp g$ if both $f = O(g)$ and $g = O(f)$.

\section{Tools from the literature}\label{sec:tools-from-liter}

The following theorem is due to Dyer--Harman~\cite[Theorem 1]{DY86}, see also Harman's monoraph~\cite[Chapter~3,
Theorem~3.2]{Harman}.

\begin{dhtheorem}\label{lem:gcdsum}
  Let  $\eps>0$. Then
  \[
    \sum_{q,r \in \mathcal{M}}\frac{\gcd(q,r)}{\sqrt{qr}}\sqrt{\psi(q)\psi(r)} \ll_{\varepsilon} |\calM|^{\varepsilon}
    \sum_{q \in \mathcal{M}}\psi(q)
  \]
  for all finite sets $\mathcal{M}\subseteq \NN$ and any $\psi:\mathcal{M}\to[0,1]$.
\end{dhtheorem}

\begin{remark*}
  In fact sharper estimates are known than the simple
  $\lvert \mathcal{M}\rvert^\varepsilon$, see Aistleitner--Berkes--Seip~\cite{ABS},
  Lewko--Radziwi\l\l~\cite{LewkoRadz},
  Bondarenko--Seip~\cite{BondarenkoSeip}, or de la
  Bretèche--Tenenbaum~\cite{BT19} for more recent developments. For
  our purposes here, the bound stated above suffices.
\end{remark*}

The following zero-one law is due to
Beresnevich--Haynes--Velani~\cite[Theorem~1]{BHV}.

\begin{bvhtheorem}\label{thm:BHV}
  Let $n\geq 1$ and $\psi:\NN\to[0,\infty)$. Then $\lambda_n(\calG_n(\psi))\in\set{0,1}$.
\end{bvhtheorem}

The following result can be found in~\cite{chungerdos}.

\begin{celemma}
  \phantomsection\label{thm:CE}
  For measurable sets $A_1, A_2, \dots, A_N\subset (X,\mu)$ in an arbitrary probability space,
  \begin{equation*}
    \mu\parens*{\bigcup_{n=1}^N A_n} \geq \frac{\parens*{\sum_{n=1}^N\mu(A_n)}^2}{\sum_{m,n=1}^N\mu(A_m\cap A_n)}.
  \end{equation*}
\end{celemma}

\section{Preliminary lemmas}\label{prelim_lem}

In this section we gather elementary calculations that go into the proof of Theorem \ref{main_thm}.
We start with the following lemma that concerns intersections of sets $S\subset\RR^n$ having the
property that
\begin{equation}\label{eq:starproperty}
  \bx=(x_1, \dots, x_n) \in S \quad\textrm{and}\quad |x_i'|\leq |x_i|\quad  (i=1, \dots, n) \quad \implies \quad  \bx'=(x_1', \dots, x_n') \in S.
\end{equation}
It is a direct adaptation of~\cite[Lemma~1]{GallagherMultiplicative}
to our setup, and it only differs from Gallagher's lemma in that our
sets are symmetric over $\bzero$, while Gallagher's lie in the positive
orthant. This is purely an expository choice; the two points of view
are equivalent.

\begin{lemma}[{\cite[Lemma~1]{GallagherMultiplicative}}]\label{lem:star_prop}
  If $A,B \subset \RR^n$ satisfy~(\ref{eq:starproperty}), then
  $\lambda(A\cap (B+\bt))$ is a decreasing function of
  $|t_i|\, (1\leq i\leq n)$ when the other components of $\bt$ are held
  fixed.
\end{lemma}

\begin{proof}
   Write $\bt = (\bs,t)$ where
  $\bs=(t_1, \dots, t_{n-1})$ and $t = t_n$. Then
  \begin{align}
    \lambda(A\cap(B+\bt)) &= \int\lambda_{1}(A^{\br}\cap(B+(\bs,t))^{\br})\, d\br  \nonumber \\
                 &= \int\lambda_{1}(A^{\br}\cap(B^{\br-\bs}+ t))\, d\br, \label{eq:lemmaslice}
  \end{align}
    where
  \begin{equation*}
    A^\br = \set{t \in\RR : (\br, t) \in A}\qquad\textrm{and}\qquad B^{\br-\bs} = \set{t \in\RR : (\br-\bs, t) \in B}.
  \end{equation*}
  By~(\ref{eq:starproperty}), the sets $A^\br$ and $B^{\br-\bs}$ are
  intervals centered at $0$. Denoting their radii $a:=a(\br)$ and
  $b:=b(\br)$, respectively, we have
  \begin{equation*}
    \lambda_{1}(A^{\br}\cap(B^{\br-\bs}+t)) =
    \begin{cases}
      \min(2a,2b) &\textrm{if } |t| < \abs{a-b} \\
      a+b-|t| &\textrm{if } |a-b|\leq |t|  < a+b \\
      0  &\textrm{if } |t| \geq a + b,
    \end{cases}
  \end{equation*}
  (interpreting $|\infty - \infty| = \infty$). This function is
  decreasing in $|t|$ so, by~(\ref{eq:lemmaslice}), the measure
  $\lambda(A\cap(B+\bt))$ is decreasing in $|t_n|$, while holding all
  other coordinates fixed. The same argument can be
    carried out for $t_i,\, i=1, \dots, n$.
\end{proof}

Next, we show that having pairwise joint coprimality is a generic condition.

\begin{lemma}\label{moebius}
Let $n \geq 2$, and let
\[\overline{\mathbb{Z}_q^n} := \{\mathbf{a} = (a_1,\ldots,a_n) \in \mathbb{Z}_q^n:
\gcd(a_i,a_j,q) = 1 \textrm{ for all } 1 \leq i < j \leq n\}.\]
Then we have $\#\overline{\mathbb{Z}_q^n} \gg_n q^n$.
\end{lemma}

\begin{proof}

Note that $q \mapsto \#\overline{\mathbb{Z}_q^n}$ is multiplicative, so we can focus on $q = p^k$, with $p$ prime. Here we have

\[\begin{split}
\#\zq &= \#\{\mathbf{a} = (a_1,\ldots,a_n) \in \mathbb{Z}_q^n: p \nmid a_i \,\forall 1 \leq i \leq n\} \\&\qquad+ \sum_{i = 1}^n \#\{\mathbf{a} = (a_1,\ldots,a_n) \in \mathbb{Z}_q^n: p \mid a_i, p \nmid a_j \,\forall j \neq i, 1 \leq j \leq n\}
\\&= q^n\left(\left(1 - \frac{1}{p}\right)^n + \frac{n}{p}\left(1 - \frac{1}{p}\right)^{n-1}\right) = q^n\left(1 - O_n\left(\frac{1}{p^2}\right)\right).
\end{split}
\]
By the multiplicative property, we can lift this to $q = \prod_{p} p^{\nu_p}$,
obtaining
\[\#\zq = q^n \prod_{p \mid q}\left(1 - O_n\left(\frac{1}{p^2}\right)\right) \gg_n q^n.\]
\end{proof}

\begin{lemma}\label{lem:Aqmeasure}
   Let
    $\psi:\NN\to[0,\infty)$ and $n\geq 2$. For all $\psi(q) \leq 2^{-n}$, we have
  \begin{equation*}
    \lambda(G_q) \asymp q^n\lambda(H_q) \asymp \lambda(A_q) \asymp \lambda(A_q^*) \asymp \psi(q)\parens*{\log \frac{1}{\psi(q)}}^{n-1},
  \end{equation*}
  with implied constants depending only on $n$. For
    $\psi(q)\geq 2^{-n}$, we have $\lambda(G_q)=\lambda(A_q)=1$.
\end{lemma}

\begin{proof}
  The last claim follows from the fact that $[-1/2, 1/2]^n \subset G_q$ when
  $\psi(q)\geq 2^{-n}$. Hence we may assume for the rest of the proof that $\psi(q) \leq 2^{-n}$.

  We prove the claim for $G_q$ by induction on $n$. When $n=2$,
  \begin{equation*}
    \lambda(G_q^{(2)}) = 4\psi(q) +  4\int_{2\psi(q)}^{1/2}\frac{\psi(q)}{y}\,dy \asymp \psi(q)\log \frac{1}{\psi(q)},
  \end{equation*}
  establishing the base case. Now let $n\geq 3$ and assume that the claim
  holds up to $n-1$. Writing $(\bx,y)\in \RR^{n-1}\times\RR$, we
have
  \begin{align*}
    \lambda(G_q) =  \int\lambda_{n-1}(G_q^y)\, dy,
  \end{align*}
  where
  \begin{align*}
    G_q^y &= \set{\bx\in\RR^{n-1} : (\bx,y)\in G_q} \\
          &= \set*{\bx\in\brackets*{-\frac{1}{2}, \frac{1}{2}}^{n-1} : |x_1|\cdots|x_{n-1}| \leq \frac{\psi(q)}{|y|}} \\
          &= G_q^{(n-1)}(\psi_y),
  \end{align*}
  with $\psi_y = \psi/|y|$. Note that if
  $|y| \leq 2^{n-1}\psi(q)$, then $G_q^y= [-1/2, 1/2]^{n-1}$. Consequently,
  \begin{align}
    \lambda(G_q) &= 2^n \psi(q) + 2 \int_{2^{n-1}\psi(q)}^{1/2}\lambda_{n-1}(G_q^{(n-1)}(\psi_y))\, dy \label{eq:forlater}\\
           &\asymp \psi(q) + \int_{2^{n-1}\psi(q)}^{1/2} \frac{\psi(q)}{y}\parens*{\log\frac{y}{\psi(q)}}^{n-2}\, dy \nonumber\\
           &\asymp \psi(q)\parens*{\log \frac{1}{\psi(q)}}^{n-1},\nonumber
  \end{align}
  which establishes the claim for $\lambda(G_q)$. Notice that $A_q = T^{-1}(G_q)$ where $T: \TT^n\to \TT^n$ is the
 mapping $\bx\mapsto q\bx\mod 1$. Since $T$ is Lebesgue
 measure-preserving, the result for $A_q$ holds.

  For $H_q$, notice that
  \begin{align}
   q\psi(q)^{-\frac{1}{n+1}} H_q\nonumber
    &= \set*{\bx\in\brackets*{-\frac{1}{2}, \frac{1}{2}}^n : |x_1|\cdots|x_n| \leq \psi(q)^{\frac{1}{n+1}}} \\
    &= G_q(\psi^{\frac{1}{n+1}}). \label{eq:HqGq_id}
  \end{align}
  It follows that 
  \begin{equation*}
    \lambda\parens*{q\psi(q)^{-\frac{1}{n+1}} H_q} \asymp \psi(q)^{\frac{1}{n+1}}\parens*{\log \frac{1}{\psi(q)^{\frac{1}{n+1}}}}^{n-1},
  \end{equation*}
  which implies that 
  \begin{equation*}
    q^n \lambda\parens*{H_q} \asymp \psi(q)\parens*{\log \frac{1}{\psi(q)}}^{n-1},
  \end{equation*}
 establishing the claim for $q^n\lambda(H_q)$.

 For $A_q^*$, recall
 that $A_q^*$ is a disjoint union of $\#\overline{\mathbb{Z}_q^n}$ 
 translates of $H_q$. Thus, by Lemma \ref{moebius},
\[
  \lambda(A_q^*)=\#\overline{\mathbb{Z}_q^n}\,\lambda(H_q)
  \asymp_n q^n\lambda(H_q),
\]
which implies the claim for $A_q^*$ and completes the proof.
\end{proof}

\section{The overlap estimate}\label{sec:overlap}

\begin{proposition}\label{prop:overlaps}
  Let $q \neq r, \psi(q),\psi(r) \leq 2^{-n}$. Then
    \begin{equation*}
      \lambda(A_q^*\cap A_r^*) \ll  \lambda(A_q^*)\lambda(A_r^*) + \gcd(q,r)\frac{\sqrt{\psi(q)\psi(r)}}{\sqrt{qr}} \min(\psi(q),\psi(r))^{1/6}.
    \end{equation*}
\end{proposition}

\begin{proof}
  Define the weight function $w = w_{q,r}^{(n)}: \mathbb{R}^n \to \mathbb{R}_{\geq 0}$ by
  \[
    w(\mathbf{x}) := \lambda\left(H_q \cap (H_r + \mathbf{x})\right).
  \]
  Then
  \[
    \lambda\left(A_q^{*} \cap A_r^{*}\right) = \sum_{\mathbf{a} \in
      \overline{\mathbb{Z}_q^n}}\sum_{\mathbf{b} \in \overline{\mathbb{Z}_r^n}}
    \lambda\left(\left(H_q + \frac{\mathbf{a}}{q}\right) \cap \left(H_r +
        \frac{\mathbf{b}}{r}\right)\right) = \sum_{\mathbf{c} \in \mathbb{Z}^n}
    \rep(\mathbf{c}) w\left(\frac{\mathbf{c}}{\lcm(q,r)}\right),
  \]
  where
  \begin{align*}
%    |\bx| &= (|x_1|,|x_2|,\ldots,|x_n|),\\
    \rep(\mathbf{c}) &= \#\left\{(\mathbf{a},\mathbf{b}) \in \zq \times \zr: \frac{\ba}{q} - \frac{\bb}{r} = \frac{\bc}{\lcm(q,r)} \right\}.% \\
    % w &= w_{q,r}: \mathbb{R}_{\geq 0}^N \to \mathbb{R}, \quad \bx \mapsto \lambda\left(H_q\cap \left(H_r + \bx\right)\right),
  \end{align*}
%  recalling that $\bx \mapsto \lambda\left(H_q\cap \left(H_r + \bx\right)\right)$ is symmetric in all coordinates.
  The Chinese Remainder Theorem gives $\rep(\bc) \leq \gcd(q,r)^n$, and
  by the definition of $\zq$, we have $\rep(\bc) = 0$ if $\bc$ has at
  least two vanishing entries. Consequently,
  \begin{align*}
    \frac{\lambda\left(A_q^{*} \cap A_r^{*}\right)}{\gcd(q,r)^n}
      &\ll \sum_{\substack{\mathbf{c} \in \mathbb{Z}^n\\ c_i > 0,\,  i = 1,\ldots,n}}w\left(\frac{\bc}{\lcm(q,r)}\right) + 
      \sum_{i = 1}^n \sum_{\substack{\mathbf{c} \in \mathbb{Z}^n\\ c_i = 0,\,  c_j > 0,\,  j \neq i}}w\left(\frac{\bc}{\lcm(q,r)}\right) \\
      &=  \sum_{\substack{\mathbf{c} \in \mathbb{Z}^n\\ c_i > 0,\,  i = 1,\ldots,n}}w\left(\frac{\bc}{\lcm(q,r)}\right) + n  \sum_{\substack{\mathbf{c} \in \mathbb{Z}^n\\ c_1 = 0, \, c_j > 0,\,  j \geq 2}}w\left(\frac{\bc}{\lcm(q,r)}\right),
  \end{align*}
  using the symmetry of $w$. Note that $H_q$ is an intersection of two sets that both satisfy \eqref{eq:starproperty}, hence $H_q$ also satisfies \eqref{eq:starproperty}. Therefore, Lemma~\ref{lem:star_prop} shows that $w$ is a decreasing function in each coordinate $|x_i|$.
  Using this property, we obtain
    
  \[
    \sum_{\substack{\mathbf{c} \in \mathbb{Z}^n\\ c_i > 0,\,  i = 1,\ldots,n}}w\left(\frac{\bc}{\lcm(q,r)}\right) \leq \lcm(q,r)^n\int_{\RR^n} w\left(\bx\right) \mathrm{d}\bx,
  \]
  as well as
  \[\sum_{\substack{\mathbf{c} \in \mathbb{Z}^{n-1}\\ c_j > 0, \,  1 \leq j \leq n-1}}w\left(\frac{(0,\bc)}{\lcm(q,r)}\right) \leq \lcm(q,r)^{n-1} 
    \int_{\RR^{n-1}} w\left(0,\bx\right) \mathrm{d}\bx.
  \]
  Note that by~\Cref{lem:Aqmeasure}, 
  \[
    \int_{\RR^n} w\left(\bx\right) \mathrm{d}\bx
    = \int_{\mathbb{R}^n}\lambda\left({H_q} \cap \left({H_r} + \bx\right)\right) \mathrm{d}\bx = \lambda(H_q)\lambda(H_r) \asymp \frac{\lambda(A_q^*)}{q^n}\frac{\lambda(A_r^*)}{r^n},
  \]
  hence
  \begin{equation*}%\label{eq:quickref}
    \lambda\left(A_q^{*} \cap A_r^{*}\right)
    \ll \lambda(A_q^*)\lambda(A_r^*) +
    q^{n-1}r^{n-1} \gcd(q,r) \int_{\RR^{n-1}} w\left(0,\bx\right) \mathrm{d}\bx.
  \end{equation*}
  The proposition follows by applying the following lemma.
\end{proof}

\begin{lemma}\label{lem:intHqHr}
  Let $n \geq 2, q,r \in \mathbb{N}, 0 \leq \psi(q),\psi(r) \leq 2^{-n}$. Then we have
  \begin{equation*}
    q^{n-1}r^{n-1} \gcd(q,r) \int_{\RR^{n-1}} w\left(0,\bx\right) \mathrm{d}\bx \ll   \gcd(q,r)\frac{\sqrt{\psi(q)\psi(r)}}{\sqrt{qr}} \min(\psi(q),\psi(r))^{1/6}.
  \end{equation*}
\end{lemma}

\begin{proof}
  Observe that with $H_q^y$ defined analogous to $G_q^y$ in the proof of Lemma
  \ref{lem:Aqmeasure},
  \begin{align}
    \int_{\RR^{n-1}} w\left(0,\bx\right) \mathrm{d}\bx
    &= \int_{\RR^{n-1}}  \lambda_n\left(H_q \cap (H_r + (0,\mathbf{x}))\right)
      \mathrm{d}\bx = \int_{\RR^{n-1}}  \int_{\mathbb{R}} \lambda_{n-1}\left(H_q^y \cap (H_r^y + \mathbf{x})\right) \mathrm{d}y\,\mathrm{d}\bx \nonumber \\
    &=  \int_{\RR} \int_{\RR^{n-1}} \lambda_{n-1}\left(H_q^y \cap (H_r^y + \mathbf{x})\right) \mathrm{d}\bx\, \mathrm{d}y = \int_{\R} \lambda_{n-1}\left(H_q^y\right)\lambda_{n-1}\left(H_r^y\right) \mathrm{d}y \nonumber\\
    &\leq \parens*{\int_{\R} \lambda_{n-1}\left(H_q^y\right)^2 \mathrm{d}y}^{1/2} \parens*{\int_{\R} \lambda_{n-1}\left(H_r^y\right)^2 \mathrm{d}y}^{1/2}, \label{eq:CS}
  \end{align}
  by the Cauchy--Schwarz inequality. Now,
  \begin{align*}
    \int_{\R} \lambda_{n-1}\left(H_q^y\right)^2 \mathrm{d}y
    &= q^{-2(n-1)}\psi(q)^{\frac{2(n-1)}{n+1}}    \int_{\R} \lambda_{n-1}\left((q\psi(q)^{-\frac{1}{n+1} }H_q)^{q\psi(q)^{-1/(n+1)}y}\right)^2 \mathrm{d}y \\
    &= q^{-2(n-1)-1}\psi(q)^{\frac{2(n-1)+1}{n+1}}    \int_{\R} \lambda_{n-1}\left((q\psi(q)^{-\frac{1}{n+1} }H_q)^{u}\right)^2 \mathrm{d}u \\
    &= q^{-2n+1}\psi(q)^{\frac{2n-1}{n+1}}    \int_{\R} \lambda_{n-1}(G_q(\psi^{\frac{1}{n+1}})^u)^2 \mathrm{d}u
  \end{align*}
 by~(\ref{eq:HqGq_id}). Continuing as in~(\ref{eq:forlater}), we bound the above by
  \begin{align*}
    &\ll q^{-2n+1}\psi(q)^{\frac{2n-1}{n+1}}\left(2^{n-1}\psi(q)^{\frac{1}{n+1}}  +   \int_{2^{n-1}\psi(q)^{\frac{1}{n+1}}}^{\infty} \parens*{\frac{\psi(q)^{1/(n+1)}}{u}\parens*{\log \frac{u}{\psi(q)^{1/(n+1)}}}^{n-2}}^2 \mathrm{d}u \right)\\
    &\ll q^{-2n+1}\psi(q)^{\frac{2n-1}{n+1}} \left(\psi(q)^{\frac{1}{n+1}} + \psi(q)^{\frac{1}{n+1}}   \int_{1}^{\infty} \parens*{\frac{1}{v}\parens*{\log v}^{n-2}}^2 \mathrm{d}v \right)\\
    &\ll q^{-2n+1}\psi(q) \psi(q)^{\frac{n-1}{n+1}} \ll q^{-2n+1}\psi(q) \psi(q)^{\frac{1}{3}}.
  \end{align*}
  Similarly, of course, 
  \begin{align*}
    \int_{\R} \lambda_{n-1}\left(H_r^y\right)^2 \mathrm{d}y \ll r^{-2n+1}\psi(r) \psi(r)^{\frac{1}{3}}.
  \end{align*}
  Putting these into~(\ref{eq:CS}) gives
  \begin{equation*}
    q^{n-1}r^{n-1} \gcd(q,r)    \int_{\RR^{n-1}} w\left(0,\bx\right) \mathrm{d}\bx \ll   \gcd(q,r)\frac{\sqrt{\psi(q)\psi(r)}}{\sqrt{qr}} \min(\psi(q),\psi(r))^{\frac{1}{6}},
  \end{equation*}
  proving the lemma.
\end{proof}

\begin{remark*} The proof shows that in fact in place of $\frac{1}{6}$, the exponent $\frac{n-1}{2n+2} = \frac{1}{2} + o(1)$ is admissible in the above. This follows the recurring
    theme in metric Diophantine approximation that higher-dimensional
    problems are easier to prove because of the smaller overlaps of
    their concomitant sets.
\end{remark*}

\section{Proof of Theorems~\ref{main_thm} and~\ref{thm:nostar_}}\label{sec:final}

\subsection*{Proof of Theorem \ref{thm:nostar_}}

\begin{proof}
  Let $S_q = [-\psi(q), \psi(q)]\times [-\frac12,\frac12]^{n-1}$, where $\psi$
  is a function such that $\sum \psi(q)=\infty$ and such that almost
  no real number is $\psi$-approximable. That is, $\psi$ is a
  Duffin--Schaeffer-type counterexample, such as was first constructed
  in~\cite{duffinschaeffer}.

  Since $S_q$ is a product of symmetric intervals, it
  satisfies~(\ref{eq:starproperty}). Notice that
  $\lambda_n(S_q) = 2\psi(q)$, hence 
  $\sum \lambda_n(S_q)$ diverges. Finally, notice that
  \begin{equation*}
    \calW = \calW_1(\psi)\times[0,1]^{n-1}, 
  \end{equation*}
  where $\calW_1(\psi)$ denotes the set of $\psi$-approximable
  numbers. Therefore $\lambda_n(\calW)=0$ by Fubini's theorem.
\end{proof}

\subsection*{Proof of Theorem~\ref{main_thm}}

The convergence part of~\Cref{main_thm} is
standard: by~\Cref{lem:Aqmeasure}, the convergence assumption in the
theorem statement implies that $\sum_{q\geq 1} \lambda(A_q)$ converges, so, by the
Borel--Cantelli lemma, the limsup set $\calG_n(\psi)$ has zero measure. For
the divergence part of Theorem~\ref{main_thm}, assume
\begin{equation*}
  \sum_{q=1}^\infty \psi(q) \parens*{\log\frac{1}{\psi(q)}}^{n-1} = \infty.
\end{equation*}
Since $\calG_n(\psi)$ contains the \emph{simultaneously}
  $\psi^{1/n}$-approximable points, the result would follow as a
  consequence of~\cite[Theorem~1]{Gallagherkt} if $\sum \psi(q)$
  diverged. We may therefore assume $\sum\psi(q)$ converges and, in
  particular, that $\psi(q)\to 0$. This implies that 
\begin{equation*}
\lim_{q\to\infty}  \psi(q) \parens*{\log\frac{1}{\psi(q)}}^{n-1}=0.
\end{equation*}
We
may now find arbitrarily large $X<Y$ such that
\begin{equation}\label{eq:sumrange}
  1 \leq   \sum_{X\leq q \leq Y} \psi(q) \parens*{\log\frac{1}{\psi(q)}}^{n-1} \leq 2, 
\end{equation}
with $\max_{X \leq q \leq Y}\psi(q) \leq 2^{-n}$ and in particular,
\begin{equation}
  \label{eq:Aqsumrange}
  \sum_{X\leq q \leq Y}\lambda(A_q^*) \asymp 1
\end{equation}
by~\Cref{lem:Aqmeasure}. We have
\begin{align*}
  \sum_{X\leq q,r\leq Y}\lambda(A_q^*\cap A_r^*)
  &\ll \sum_{X\leq q\leq Y}\lambda(A_q^*) + \sum_{X\leq r < q \leq Y} \lambda(A_q^*\cap A_r^*)  \\
  &\ll \parens*{\sum_{X \leq q\leq Y}\lambda(A_q^*)}^2 + \sum_{X\leq r < q \leq Y} \frac{\gcd(q,r)}{\sqrt{qr}}\sqrt{\psi(q)\psi(r)}\min(\psi(q),\psi(r))^{1/6},
\end{align*}
 by~\Cref{prop:overlaps}. Now
\begin{multline*}
  \sum_{X\leq r < q \leq Y} \frac{\gcd(q,r)}{\sqrt{qr}}\sqrt{\psi(q)\psi(r)}\min(\psi(q),\psi(r))^{1/6} \\
  \leq  \sum_{j\geq 1}  \sum_{\substack{X\leq r < q \leq Y \\ 2^{-j}\leq \min(\psi(q),\psi(r))\leq 2^{-j+1}}} \frac{\gcd(q,r)}{\sqrt{qr}}\sqrt{\psi(q)\psi(r)}\min(\psi(q),\psi(r))^{1/6} \\
  \ll  \sum_{j\geq 1} 2^{-j/6} \sum_{\substack{X \leq r < q \leq Y \\ \psi(q),\psi(r)\geq 2^{-j}}} \frac{\gcd(q,r)}{\sqrt{qr}}\sqrt{\psi(q)\psi(r)} \\
    \ll  \sum_{j\geq 1} 2^{-j/6}\parens*{\#\set*{q \in [X,Y]: \psi(q) \geq 2^{-j}}}^{\eps} \sum_{q\in [X,Y]} \psi(q)
\end{multline*}
by~\Cref{lem:gcdsum} (\Cref{sec:tools-from-liter}).
From~(\ref{eq:sumrange}), we obtain
\begin{equation*}
  \#\set*{q \in [X,Y]: \psi(q) \geq 2^{-j}} \ll 2^j, 
\end{equation*}
therefore, for $\varepsilon < \tfrac{1}{6}$,
\begin{align*}
  \sum_{X\leq q,r\leq Y}\lambda(A_q^*\cap A_r^*)
  &\ll  \parens*{\sum_{X \leq q\leq Y}\lambda(A_q^*)}^2 + \parens*{ \sum_{j\geq 1} 2^{-j(1/6 - \varepsilon)}} \sum_{q\in [X,Y]} \psi(q) \\
  &\ll \parens*{\sum_{X \leq q\leq Y}\lambda(A_q^*)}^2 +  \sum_{q\in [X,Y]} \psi(q) \ll \parens*{\sum_{X \leq q\leq Y}\lambda(A_q^*)}^2.
\end{align*}
It follows by the~\hyperref[thm:CE]{Chung--Erd{\H o}s Lemma}
(\Cref{sec:tools-from-liter}) that
\begin{equation*}
  \lambda\parens*{\bigcup_{X\leq q\leq Y} A_q^*} \gg 1,
\end{equation*}
and since there are arbitrarily large $X < Y$
satisfying~(\ref{eq:Aqsumrange}), this shows that
\begin{equation*}
  \lambda\parens*{\limsup_{q\to\infty} A_q^*} >0,
\end{equation*}
hence $\lambda(\calG_n(\psi))>0$. By~\Cref{thm:BHV} (\Cref{sec:tools-from-liter}), this implies $\lambda(\calG_n(\psi)) = 1$. \qed

\subsection*{AI declaration}

  We used ChatGPT-6 Astra to generate TikZ code
  for~\Cref{fig:AqHq,fig:Aqprime,fig:Aqtilde,fig:errorcontrib}, which
  we then adjusted. We did not use AI for any of the mathematical work
  or the writing.

\bibliographystyle{plain}

\end{document}